\documentclass{amsart}

\usepackage[centertags]{amsmath} 
\usepackage{amsfonts,amssymb,amsthm,amscd,latexsym,euscript,bbm,stmaryrd}
\usepackage[all]{xy}

\newtheorem {theorem}{Theorem}[section]
\newtheorem {corollary}[theorem]{Corollary}

\newtheorem {lemma}[theorem]{Lemma}
\theoremstyle{definition}
\newtheorem {definition}[theorem]{Definition}

\newtheorem {remark}[theorem]{Remark}
\newtheorem {conjecture}[theorem]{Conjecture}

\numberwithin{equation}{section}        

\newcounter{zahl}%
    {\end{list}}%
\newcommand{\co}{\colon\thinspace}
\newcommand{\mc}[1]{\ensuremath{\mathcal{#1}}}
\newcommand{\colim}{\operatorname*{colim}}
\newcommand{\holim}{\operatorname*{holim\,}}
\newcommand{\free}{\,\_\!\_\,}
\newcommand{\id}{\operatorname{id}}
\newcommand{\ul}[1]{\underline{#1}}
\newcommand{\ol}[1]{\overline{#1}}
\newcommand{\wt}[1]{\widetilde{#1}}
\newcommand{\sgr}{\ensuremath{s{\rm Gr}}}%
\newcommand{\Ev}{\ensuremath{{\rm Ev}}}%
\newcommand{\Lie}{\ensuremath{{\rm Lie}}}%
\newcommand{\Nil}{\ensuremath{{\rm Nil}}}%
\newcommand{\Alg}{\ensuremath{{\rm Alg}}}%
\newcommand{\Sfinp}{\ensuremath{\mc{S}_*^{{\rm fin}}}}%
\newcommand{\map}{\ensuremath{{\rm map}}}
\newcommand{\dgrm}[1]{\ensuremath{\smash{\underset{\widetilde{\hphantom{#1}}}{#1}} \mathstrut}}%
\newcommand{\diagr}[1]{ \begin{equation*} \xymatrix{#1} \end{equation*}}%
 
\usepackage{hyperref} 

\begin{document}

\title  [Homotopy nilpotent groups]{Homotopy nilpotent groups}
\author{Georg Biedermann}
\author{William G. Dwyer}


\address{Max-Planck-Institut f\"ur Mathematik, Vivatsgasse 7, 53111 Bonn, Germany}
\address{Department of Mathematcs, 255 Hurley, University of Notre Dame, Notre Dame IN 46556 USA}

\email{gbiederm@mpim-bonn.mpg.de}
\email{dwyer.1@nd.edu}



\date{\today}

\begin{abstract}
We study the connection between the Goodwillie tower of the identity and the lower central series of the loop group on connected spaces. We define the simplicial theory of homotopy $n$-nilpotent groups. This notion interpolates between infinite loop spaces and loop spaces.
We prove that the set-valued algebraic theory obtained by applying $\pi_0$ is the theory of ordinary $n$-nilpotent groups and that the Goodwillie tower of a connected space is determined by a certain homotopy left Kan extension. We prove that $n$-excisive functors of the form $\Omega F$ have values in homotopy $n$-nilpotent groups. 
\end{abstract}

\maketitle

\section{Introduction}

The aim of this article is to explore the connection between the Goodwillie tower of the identity functor and the lower central series of Kan's loop group of a connected space. We express it with the aid of simplicial algebraic theories. We expect the reader to be familiar with the basic notions of homotopical algebra and Goodwillie's calculus of homotopy functors. We define the notion of homotopy $n$-nilpotent groups.
The main theorems explain their relation to loop spaces \ref{loop}, to infinite loop spaces \ref{infinite loop}, to ordinary nilpotent groups \ref{pi0-iso}, to the Goodwillie tower of the identity \ref{main theorem}, and tell us that $n$-excisive functors of the form $\Omega F$ take values in the category of homotopy $n$-nilpotent groups \ref{values}. 

Let us introduce some notation valid for the rest of the article. Let $\mc{S}_*$ be the category of pointed simplicial sets. 
Let $\mc{S}_0$ be the category of reduced simplicial sets, i.e. simplicial sets with exactly one $0$-simplex. Let $X$ be an object in $\mc{S}_0$. Further let \mc{F} denote the $\mc{S}_*$-category of $\mc{S}_*$-functors from finite pointed simplicial sets to $\mc{S}_*$. A homotopy functor in \mc{F} is a functor that preserves weak equivalences.

For a homotopy functor $\dgrm{X}$ in \mc{F} Goodwillie \cite{Goo:calc3} constructs a tower of functors
    $$ \dgrm{X} \to ... \to P_n\dgrm{X}\to P_{n-1}\dgrm{X}\to ...\to P_1\dgrm{X}\to P_0\dgrm{X}=\dgrm{X}(\ast), $$
where the $n$-th stage is the universal $n$-excisive homotopy functor under \dgrm{X}. Here, $n$-excision is a higher version of excision; a $1$-excisive functor is a homotopy functor with a Mayer-Vietoris sequence. For $\dgrm{X}=\id$ the Goodwillie tower converges on simply connected spaces to the identity \cite{Goo:calc2}, but on connected spaces it converges to the Bousfield-Kan completion \cite{Arone-Kankaan:iterated-Snaith}:
    $$ \holim_{n} P_n(\id)(X)\simeq \mathbbm{Z}_{\infty}X $$
We have $P_1(\id)\simeq\Omega^{\infty}\Sigma^{\infty}$, the stable homotopy functor. So the Goodwillie tower interpolates between stable and unstable homotopy. The map $X\to P_n(\id)(X)$ is roughly $(n+1)c$-connected, if $X$ is $c$-connected.

Let $\sgr$ denote the category of simplicial groups. Kan's loop group functor $G$ is  part of a Quillen equivalence
    $$ G\co\mc{S}_0\leftrightarrows\sgr :\!\ol{W}. $$
The homotopy category of $\mc{S}_0$ and of connected spaces are equivalent, see \cite[V]{GoJar:simp}. 
The lower central series filtration of $GX$ was studied by Curtis \cite{Curtis:relations} who proved that for simply connected spaces the connectivity of the map
    $$ GX\to GX/\Gamma_{n+1}GX $$
increases logarithmically with $n$. The tower $\{GX/\Gamma_{n+1}GX\}_{n\ge 1}$ associated to the filtration converges to the identity on simply connected spaces. As noted by Kan:
    $$ \pi_{s-1}(GX/[GX,GX])\cong H_sX \text{ for all } s\ge 1$$
Since the $n$-th stage of the lower central series tower is $n$-excisive there is a map from the looped version $\{\Omega P_n(\id)(X)\}_{n\ge 1}$ of the Goodwillie tower of the identity at $X$ to the lower central series tower of $GX$. 

Algebraic theories were introduced by Lawvere \cite{Lawvere:semantics} to obtain categorical descriptions of algebraic structures like groups, rings, Lie algebras, etc. An algebraic theory is a category $T$ having the natural numbers $k\ge 0$ as objects such that $k$ is the product in $T$ of $k$ copies of $1$. The maps from $k$ to $1$ are to be thought of as the $k$-ary operations of $T$. They can be canonically identified with the free objects on $k$ generators.
Algebras over $T$ are product preserving functors from $T$ to sets.

For purposes in homotopy theory we need to consider simplicial algebraic theories where $T$-algebras have values in simplicial sets. These were first considered by Reedy \cite{Reedy:theories} and more recently by \cite{Schwede:theories} and \cite{Badzioch:theories}. 
It is convenient to study pointed versions where algebras are functors from $T$ to pointed simplicial sets and the category $T$ itself will be enriched over $\mc{S}_*$. If the theory has only one constant -- as in our case -- there is no loss in generality.
We also need a weaker notion of algebra: homotopy $T$-algebras were introduced by Badzioch \cite{Badzioch:theories}. They are functors from $T$ to $\mc{S}_*$ that commute with products up to homotopy. 
The free objects on $k$ generators of these simplicial theories are given by
    $$ \Omega(\bigvee_{i=1}^kS^1)\simeq\Omega\Sigma\bigvee_kS^0\, \text{ and }\, \Omega P_1(\id)(\bigvee_{i=1}^kS^1)\simeq\Omega^{\infty}\Sigma^{\infty}\bigvee_kS^0. $$ 
We define for $1\le n<\infty$ new theories $\mc{P}_n$ with free objects given by
    $$ \Omega P_n(\id)(\bigvee_{i=1}^kS^1). $$
Homotopy $n$-nilpotent groups are defined as homotopy $\mc{P}_n$-algebras.
For $n=\infty$ we get back loop spaces \ref{loop} and for $n=1$ infinite loop spaces \ref{infinite loop}.
Our next result says that by applying $\pi_0$ we get back the ordinary theory of $n$-nilpotent groups \ref{pi0-iso}.

We obtain morphisms of theories
   $$ \mc{P}_{\infty}\to ...\to\mc{P}_n\to\mc{P}_{n-1}\to ...$$
induced by the maps in the Goodwillie tower and pullback functors
   $$ \varphi_n\co\mc{S}_*^{\mc{P}_n}\to\mc{S}_*^{\mc{P}_{\infty}}, $$
with left adjoints
   $$ \lambda_n\co\mc{S}_*^{\mc{P}_{\infty}}\to\mc{S}_*^{\mc{P}_n}. $$

For a connected space $X$ we have $\Omega X\simeq\dgrm{X}(1^+)$ where $\dgrm{X}\co\mc{P}_\infty\to\mc{S}_*$ is a homotopy $\infty$-algebra. 
Our second result \ref{main theorem} exhibits the Goodwillie tower of the identity $\Omega P_n(\id)(X)$ as the homotopy left Kan extension along $\varphi_n$:
   $$ \Omega P_n(\id)(X)\simeq (L\lambda_n\dgrm{X})(1^+). $$

Finally we prove in \ref{values} that functors of the form $\Omega F$ where $F$ is $n$-excisive, naturally take vakues in the category of homotopy $n$-nilpotent groups. It follows \ref{enriched} that functors of this form are naturally enriched in homotopy $n$-nilpotent groups.
This justifies a closer study of the notion.

{\bf Acknowledgments:} We would like to thank Andr\'{e} Joyal and Gerald Gaudens for several stimulating discussions and Bernard Badzioch for explaining us theorem \ref{badzioch}. The research for this paper was started while both authors were guests at the Thematic Program on Geometric Applications of Homotopy Theory at the Fields Institute for Mathematics, Toronto. The first author was a post-doctoral fellow at the University of Western Ontario, London, Ontario. The paper was finished while the first author was guest at the Max-Planck-Institut f\"ur Mathematik in Bonn.

\section{Simplicial algebraic theories}
\label{section:theories}

We will consider theories enriched over the category $\mc{S}_*$ of pointed simplicial sets and algebras with values in $\mc{S}_*$. As it turns out in our case \ref{pointed vs unpointed}, the resulting notions are equivalent. Let us recast the definitions in the pointed version.
For a simplicial category \mc{C} we will refer to its simplicial set of morphisms by $\mc{C}(\free,\free)$.

\begin{definition}\label{Gamma}
Let $\Gamma$ be the opposite of the category of finite pointed sets.
The category $\Gamma$ has all products and every object is isomorphic to an object of the form
    $$ k^+=\{1,...,k\}\cup\{+\}.$$
For every $1\le s\le k$ we have maps $i_s^k\co k^+\to 1^+$ given by the inclusion of the pointed set $1^+$ to $k^+$ where the non basepoint of $1^+$ maps to $s\in k^+$.
These maps induce an isomorphism
    $$ \prod_{s=1}^ki_s^k\co k^+\cong \prod_{s=1}^k 1^+.$$
We can view $\Gamma$ as a discrete simplicial category.
\end{definition}

\begin{definition}
A {\it simplicial pointed algebraic theory} is a category $T$ enriched over pointed simplicial sets $\mc{S}_*$ having the same discrete set of objects as $\Gamma$ together with a functor $\Gamma\to T$, which is the identity on objects and preserves products. We will abbreviate this as {\it simplicial theory}. Morphisms of simplicial theories are product preserving $\mc{S}_*$-functors under $\Gamma$. 
\end{definition}

The category $T$ is usually given as a full subcategory of some other category. So the morphisms in $T$ are often left understood, and we will often confuse the objects $k^+$ of $T$ with their images under this full inclusion. If we want to emphasize a particular theory $T$, we will denote the objects by $T(k^+)$.

\begin{definition}
A {\it strict $T$-algebra} is a simplicial functor $\dgrm{X}\co T\to\mc{S}_*$ that preserves products strictly. This means that the map 
    $$ \prod_{s=1}^k\dgrm{X}(i_s^k)\co\dgrm{X}(k^+)\to\dgrm{X}(1^+)^k $$
is an isomorphism.
A {\it homotopy $T$-algebra} is a simplicial functor from $T$ to $\mc{S}_*$ that preserves products up to weak equivalence. This means that the above map is a weak equivalence. 
\end{definition}

The category of strict $T$-algebras is a reflexive subcategory of the category of $\mc{S}_*$-functors $\mc{S}_*^T$ from $T$ to $\mc{S}_*$ and carries a model structure where fibrations and weak equivalences are detected by the forgetful functor to $\mc{S}_*^T$. This is explained in the unpointed case in \cite{Badzioch:theories}.
We will examine, how to relate these situations.

\begin{definition}
The {\it constants} of a simplicial theory are the $0$-ary operations, i.e. the simplicial set
   $$ A_0=T(0^+,1^+). $$
A {\it simplicial theory with one constant} is a simplicial theory such that $A_0\cong\ast$. 
\end{definition}

An example of an ordinary algebraic theory that has more than one constant, is the theory of rings with $A_0=\mathbbm{Z}$. However, the theory of groups and $n$-nilpotent groups for $n\ge 1$ has only one constant.

\begin{remark}
Let $T$ be a simplicial algebraic theory. Then the forgetful functor $u\co\mc{S}_*\to\mc{S}$ induces a functor
    $$ u_*\co\mc{S}_*^T\to\mc{S}^T. $$
This restricts to a functor 
    $$ u_*\co\Alg_{T,*}\to\Alg_T $$
from the category $\Alg_{T,*}$ of pointed $T$-algebras to $T$-algebras $\Alg_T$. 
\end{remark}

\begin{lemma}\label{one constant}
If $T$ is a simplicial theory with one constant, the functor
    $$ u_*\co\Alg_{T,*}\to\Alg_T $$
is an isomorphism of model categories.
\end{lemma}

\begin{proof}
Given an unpointed $T$-algebra $X$, we can always supply it with a canonical basepoint
    $$ \ast\cong X(0^+)\to X(1^+),$$
induced by the unique constant $0^+\to 1^+$ in $T$. We obtain an inverse functor for $u_*$. Because weak equivalences and fibrations are given in both model categories by the ones on underlying simplicial sets, we have an isomorphism of model categories.
\end{proof}

Now we want to consider homotopy $T$-algebras.
The category $\mc{S}_*^{T}$ can be equipped with a model structure where the objectwise fibrant homotopy $T$-algebras are exactly the fibrant objects. This model structure is a localization of the projective model structure on $\mc{S}_*^T$.

\begin{definition}
We will call the category $\mc{S}_*^T$ together with this localized model structure the {\it homotopy algebra model structure} $\left(\mc{S}_*^T\right)_{\rm halg}$.
\end{definition}

It is shown in \cite{Badzioch:theories} that there is a Quillen equivalence between strict $T$-algebras and homotopy $T$-algebras. This striking result tells us that -- independently of the theory $T$ -- any homotopy $T$-algebra can be rigidified. 

\begin{theorem}[Badzioch]\label{badzioch}
Let $F\co S\to T$ be a morphism of simplicial theories. If $F$ is a weak equivalence of simplicial categories then pulling back along $F$ is the right adjoint of a Quillen euqivalence between the associated homotopy algebra model categories.
\end{theorem}

\begin{proof}
In theorem 2.1 of \cite{DK:equivalence} it is shown that $F^*\co\mc{S}_*^T\to\mc{S}_*^S$ is the right adjoint of a Quillen equivalence between the projective model structures. The homotopy algebra model structures are left Bousfield localizations. This process preserves Quillen equivalences by theorem 3.3.20 from \cite{Hir:loc}.
\end{proof}

\section{The $n$-excisive model structure}

\begin{definition}\label{mc{F}}
We denote by \mc{F} the category of $\mc{S}_*$-enriched functors from finite pointed simplicial sets $\Sfinp$ to pointed simplicial sets $\mc{S}_*$.
\end{definition}

The category \mc{F} is enriched, tensored and cotensored over $\mc{S}_*$ where both tensor and cotensor are given objectwise. It carries a projective model structure where weak equivalences and fibrations are given objectwise. 

For an introduction to Goodwillie's calculus of homotopy functors and in particular for the notion of $n$-excisive homotopy functor we refer to \cite{Goo:calc3} and \cite{Kuhn:overview}.

In \cite{BCR:calc} and \cite{Dwyer:localizations} the projective model structure on \mc{F} was localized to obtain the $n$-excisive model structure where the fibrant objects are exactly the $n$-excisive homotopy functors. 
A map $\dgrm{X}\to \dgrm{Y}$ is an $n$-excisive weak equivalence if and only if it induces an objectwise weak equivalence
    $$ P_n\dgrm{X}\to P_n\dgrm{Y}. $$
Here $P_n\dgrm{X}$ denotes the $n$-th stage in the Goodwillie tower of the functor $\dgrm{X}^h$ which is the functor \dgrm{X} pre- and postcomposed with an objectwise fibrant replacement functor in $\mc{S}_*$. 
However it is more convenient for us to consider the injective model structure on \mc{F} constructed by Joyal \cite{Joyal:letter} and Jardine \cite{jar:simplicial-presheaves} where cofibrations are given by all inclusions. This model structure is also proper and simplicial with the advantage that all objects are cofibrant. The same techniques as in \cite{BCR:calc} apply to arrive at an $n$-excisive model structure on \mc{F} where a map $\dgrm{X}\to \dgrm{Y}$ is an $n$-excisive equivalence as above and an an $n$-excisive fibration if and only if it is an injective fibration such that the square
\diagr{ \dgrm{X} \ar[r]\ar[d] & P_n\dgrm{X} \ar[d] \\ \dgrm{Y} \ar[r] & P_n\dgrm{Y}}
is an objectwise homotopy pullback square. We have:

\begin{theorem}\label{inj. n-exc. model str.}
The injective $n$-excisive model structure on \mc{F} is a cofibrantly generated proper simplicial model structure.
\end{theorem}

\begin{corollary}\label{n-excisive and mapping space}
Let $\dgrm{Y}$ be an injectively fibrant $n$-excisive homotopy functor. Then for every $\dgrm{X}$ in \mc{F} we have a natural weak equivalence
    $$ \mc{F}(P_n\dgrm{X},\dgrm{Y})\simeq\mc{F}(\dgrm{X},\dgrm{Y}). $$
\end{corollary}

\section{Homotopy $n$-nilpotent groups}

We will now describe the simplicial theory of homotopy $n$-nilpotent groups. 

\begin{definition}
In the category \mc{F} let $(\free)^{\rm inj}$ be a fibrant replacement functor with respect to the injective model structure.
\end{definition}

\begin{definition}
We define a full subcategory $\mc{P}_n$ of the category $\mc{F}$, which has for each natural number $k\ge 0$ exactly one object given by
    $$ \mc{P}_n(k^+)=\prod_{i=1}^k\Omega(P_n(\id))^{\rm inj}. $$
We also define for $n=\infty$ the category $\mc{P}_{\infty}$ with objects given by
    $$ \mc{P}_\infty(k^+)=\prod_{i=1}^k\Omega(\id)^{\rm inj}. $$
We employ the convention that the empty product is the final object $\ast$, and so we have for all $1\le n\le\infty$ and $0\le k<\infty$:
    $$ \mc{P}_n(k^+)\cong\prod_{i=1}^k\mc{P}_n(1^+).$$
We let $I_n\co\mc{P}_n\to\mc{F}$ be the inclusion functor. 

We also define a functor $\gamma_n\co\Gamma\to\mc{P}_n$ on objects simply by $k^+\mapsto\mc{P}_n(k^+)$. For a morphism $f^{\rm op}\co\ell^+\to k^+$ in $\Gamma$ which is represented by a map $f\co k^+\to\ell^+$ of pointed sets we define $\gamma_n(f^{\rm op})$ on the $i$-th factor of $\mc{P}_n(\ell^+)$ as
\diagr{ \mc{P}_n(k^+)=\prod_k\Omega (P_n(\id))^{\rm inj}\ar[rr]_-{{\rm pr}_{f(i)}} && \Omega (P_n(\id))^{\rm inj} }
the projection on the $f(i)$-th factor. 

The natural transformation $\id\to P_n(\id)$ induces a morphism of theories 
    $$ p_n\co\mc{P}_{\infty}\to\mc{P}_n .$$
Obviously we have the equation $p_n\gamma_{\infty}=\gamma_n$.
\end{definition}

\begin{remark}
As a full subcategory of \mc{F} the category $\mc{P}_n$ for $1\le n\le \infty$ is en\-riched over $\mc{S}_*$.
The category $\mc{P}_n$ constitutes a simplicial theory as discussed in section \ref{section:theories}. 
Therefore we can con\-sider $\mc{P}_n$-algebras and homotopy $\mc{P}_n$-algebras. Objectwise fibrant homotopy $\mc{P}_n$-algebras are the fibrant objects in $(\mc{S}_*^{\mc{P}_n})_{\rm halg}$.
\end{remark}

\begin{definition}
A {\it homotopy $n$-nilpotent group} is a space weakly equivalent to $\dgrm{X}(1)$, for some homotopy $\mc{P}_n$-algebra \dgrm{X}. 
\end{definition}

\begin{lemma}\label{morphisms in P_n}
For all $n\ge 1$ and $k,\ell\ge 0$ we have canonical weak equivalences
    $$ \mc{P}_n(k^+,\ell^+)\simeq \left(\Omega P_n(\id)(\bigvee_{i=1}^k S^1)\right)^\ell$$
\end{lemma}

\begin{proof}
Recall that \mc{F} from \ref{mc{F}} is endowed with the injective $n$-excisive model structure where all objects are cofibrant.
We compute:
\begin{align*}
    \mc{P}_n(k^+,\ell^+)  &\cong\mc{F}(\prod_k\Omega (P_n(\id))^{\rm inj},\prod_{\ell}\Omega(P_n(\id))^{\rm inj}) \\
                       &\simeq\mc{F}(\Omega(\id)^k,\Omega (P_n(\id))^{\rm inj})^\ell \\
                       &\cong\mc{F}(\map_{\mc{S}_*}(\bigvee_k S^1,\free),\Omega (P_n(\id))^{\rm inj})^\ell \\
                       &\cong\left(\Omega (P_n(\id))^{\rm inj}(\bigvee_k S^1)\right)^\ell \\
                       &\simeq\left(\Omega P_n(\id)(\bigvee_k S^1)\right)^\ell
\end{align*}
The weak equivalence in step 2 comes from \ref{n-excisive and mapping space}.
We also use the enriched Yoneda lemma.
\end{proof}

\begin{remark}\label{pointed vs unpointed}
By \ref{morphisms in P_n} the theory $\mc{P}_n$ has only one constant.
So by lemma \ref{one constant} we can work in the pointed setting without losing information.
\end{remark}

\begin{corollary}\label{free objects}
The free homotopy $n$-nilpotent group on $k$ generators is given by 
    $$\Omega P_n(\id)(\bigvee_{i=1}^kS^1).$$
\end{corollary}

\begin{proof}
The free algebra $A_k$ on $k$ generators in any simplicial theory $T$ can be obtained by the following formula:
   $$ A_k\cong T(k^+,1^+) .$$
Now the statement follows from \ref{morphisms in P_n}.
\end{proof}

If $F_k$ is the free group on $k$ generators then we have the following canonical weak equivalences:
\begin{equation}\label{equivalence}   
    \mc{P}_\infty(k^+,1^+)\simeq\Omega\bigvee_kS^1\simeq\Omega BF_k\simeq F_k 
\end{equation}
It follows that the theory $\mc{P}_{\infty}$ is weakly equivalent as a simplicial category to the discrete theory of groups.

\begin{theorem}\label{loop}
The category $(\mc{S}_*^{\mc{P}_{\infty}})_{\rm halg}$ is Quillen equivalent to the theory of simplicial groups.
In particular, the homotopy category of homotopy $\infty$-nilpotent groups is equivalent to the homotopy category of loop spaces.
\end{theorem}

\begin{proof}
Follows readily from the equivalences (\ref{equivalence}) and theorem \ref{badzioch}.
\end{proof}

In \cite{Badzioch-Chung-Voronov} the authors define a theory $\mc{T}_n$ such that its homotopy algebras are exactly $n$-fold loop spaces.

\begin{definition}
Set $\mc{T}_n(k^+)=\bigvee_kS^n$ and take as morphisms the derived mapping space
    $$ \mc{T}_n(k^+,\ell^+)=\map_{\mc{S}_*}^{\rm der}(\bigvee_{\ell}S^n,\bigvee_kS^n). $$
So $\mc{T}_n$ is the opposite of the full subcategory of $\mc{S}_*$ given by the finite wedges of $n$-spheres. 
The suspension functor induces a morphism $t_n\co\mc{T}_n\to\mc{T}_{n+1}$ of simplicial theories:
    $$ \mc{T}_n(k^+,\ell^+)=(\Omega^n\bigvee_kS^n)^{\ell}\to(\Omega^{n+1}\bigvee_kS^{n+1})^{\ell}\cong\mc{T}_{n+1}(k^+,\ell^+). $$
\end{definition}

\begin{theorem}[Thm 1.1, \cite{Badzioch-Chung-Voronov}]\label{Badzioch-Chung-Voronov}
A pointed space $X$ is an $n$-fold loop space if and only if there exists a homotopy $\mc{T}_n$-algebra $\dgrm{X}$ with $\dgrm{X}(1^+)\simeq X$.
\end{theorem}

\begin{definition}
We obtain a morphism of theories $\vartheta_n\co\mc{T}_n\to\mc{P}_1$ induced by the maps
    $$ \mc{T}_n(k^+,\ell^+)=(\Omega^n\bigvee_kS^n)^{\ell}\to(\colim_{s}\Omega^{n+s}\bigvee_kS^{n+s})^{\ell}\simeq\mc{P}_1(k^+,\ell^+) $$
together with the equation $\vartheta_{n+1}t_n=\vartheta_n$.
\end{definition}

\begin{remark}\label{colimit}
The theory $\mc{P}_1$ is canonically weakly equivalent to the colimit of the sequence
\diagr{\dots\ar[r] & \mc{T}_n \ar[r]^-{t_n} & \mc{T}_{n+1} \ar[r]^-{t_{n+1}} & \dots } 
in the category of simplicial categories (with fixed set of objects $\mathbbm{N}$) and hence in the category of simplicial theories.
\end{remark}

\begin{theorem}\label{infinite loop}
A pointed space $X$ is an infinite loop space if and only if there exists a homotopy $\mc{P}_1$-algebra $\dgrm{X}$ with $\dgrm{X}(1^+)\simeq X$.
The homotopy category of homotopy $1$-nilpotent groups is equivalent to the homotopy category of infinite loop spaces.
\end{theorem}

\begin{proof}
Recall theorem \ref{Badzioch-Chung-Voronov}.
By \ref{badzioch} and remark \ref{colimit} a space $X$ is an infinite loop space if and only if the associated functor
    $$ \dgrm{X}\co\Gamma\to\mc{S}_*, k^+\mapsto X^k $$
extends via the maps $t_n$ to a product-perserving functor from $\mc{T}_n$ for all $n\ge 0$.
The existence of the morphism $\vartheta_n$ shows that $X$ is a homotopy $\mc{P}_1$-algebra if and only $\dgrm{X}$ restricts to a homotopy $\mc{T}_n$-algebra for each $n\ge 0$. The equivalence of homotopy categories now also follows.
\end{proof}

\section{The lower central series of the loop group}

\begin{definition}
Let $G$ be a group. For subgroups $H$ and $K$ of $G$ let $[H,K]$ denote the normal subgroup generated by elements of the form $h^{-1}k^{-1}hk$ where $h\in H$ and $k\in K$.
The {\it lower central filtration} for $G$ is defined in the following inductive way: Let
    $$ \Gamma_1G=G \text{ and } \Gamma_{n+1}G=[G,\Gamma_nG]. $$
We obtain a filtration of $G$ by normal subgroups with an associated tower:
\diagr{ G/\Gamma_2G \ar@{=}[d] & G/\Gamma_3G \ar[l] & G/\Gamma_4G \ar[l] & ...\ar[l]  \\
        G/[G,G] & \Gamma_2G/\Gamma_3G \ar[u] & \Gamma_3G/\Gamma_4G \ar[u] & & }
This is the {\it lower central series} of $G$. A group $G$ is called {\it $n$-nilpotent} if $\Gamma_{n+1}G=0$.
\end{definition}

\begin{definition}
For an abelian group $A$ let 
    $$\Lie_*A=\bigoplus_{n\ge 1}\Lie_nA$$ 
be the free graded Lie algebra on $A$. 
\end{definition}

\begin{remark}\label{Witt}
The Poincar\'{e}-Birkhoff-Witt theorem \cite[I.4.3]{Serre:Lie} says that there is a natural isomorphism of abelian groups
    $$ \Lie_n(G/[G,G])\cong\Gamma_nG/\Gamma_{n+1}G $$
for every free group $G$. In fact, this group is free abelian on generators given by a Hall basis of basic commutators of weight $n$ over the generators of $G$ \cite{Hall:higher-commutators}.
\end{remark}

Remember that $\mc{S}_0$ denotes the category of reduced simplicial sets and $\sgr$ the category of simplicial groups. Let $G\co\mc{S}_0\to\sgr$ be Kan's loop group functor. We can apply the functors from the lower central series degreewise. 
\begin{definition}
Let $\Gamma^n\co\mc{S}_0\to\mc{S}_*$ be the functor given by
    $$ \Gamma^nX= B\left( GX/\Gamma_{n+1}GX\right). $$
The functor $\ol{\Gamma}^n\co\mc{S}_0\to\mc{S}_*$ will be given by  
    $$ \ol{\Gamma}^nX=B\left( \Gamma_{n}GX/\Gamma_{n+1}GX\right). $$
\end{definition}

\begin{remark}
The loop group is a free simplicial group. It follows from a theorem by Dold \cite{Dold:free} that both functors $\Gamma^n$ and $\ol{\Gamma}^n$ preserve weak equivalences. Moreover with \ref{Witt} we have a formula:
    $$ \ol{\Gamma}^nX\cong\Lie_n(\wt{\mathbbm{Z}}X)$$
Here $\wt{\mathbbm{Z}}X=\mathbbm{Z}X/\mathbbm{Z}\ast$ is the reduced free simplicial abelian group on $X$.
\end{remark}

\begin{remark}\label{Curtis result}
It is proved by Curtis \cite{Curtis:relations} that for a simply connected space $X$ the map
    $$ GX\to GX/\Gamma_nGX=\Gamma^{n-1}X $$
is $\{\log_2n\}$-connected where $\{a\}$ is the least integer $\ge a$. If $X$ is merely connected, the tower $\{\Gamma^{n}(X)\}_{n\in\mathbbm{N}}$ converges to the Bousfield-Kan completion $\mathbbm{Z}_{\infty}(X)$. Compare \ref{convergence of Goodwillie tower} about the Goodwillie tower of the identity.
\end{remark}

For $n=1$ we have $\ol{\Gamma}^1X=B\wt{\mathbbm{Z}}X$. This functor is linear, because we have for all $s\ge 0$:
    $$\pi_sB\wt{\mathbbm{Z}}X\cong \wt{H}_sX$$
Here $\wt{H}_*X$ is the reduced singular homology of the reduced space $X$.
More generally, there is the following lemma.
\begin{lemma}\label{n-excisive 1}
The functor $\ol{\Gamma}^n$ is $n$-excisive.
\end{lemma}

\begin{proof}
Consider for free abelian groups $A_1,...,A_{n+1}$ the cubical diagram
    $$ P(\ul{n+1})\to{\rm FrAb}, S\mapsto \bigoplus_{i\in\ul{n+1}-S}A_i ,$$
where the maps are induced by collapsing summands.
The $(n+1)$-st cross effect of the functor $\Lie_n\co{\rm FrAb}\to{\rm Ab}$ is given by $H_0$ of the associated complex
    $$ k\mapsto \bigoplus_{|S|=k}\Lie_n(\bigoplus_{i\in\ul{n+1}-S}A_i)=:L_k. $$
But the map $L_0\to L_1$ is clearly injective.
It follows that the composition
    $$ \ol{\Gamma}^nX=B\Lie_n(\wt{\mathbbm{Z}}X) $$
is $n$-excisive.
\end{proof}

\begin{remark}
It is not true though that the functor $\ol{\Gamma}^n$ is $n$-homogeneous. By Curtis' result \ref{Curtis result} the tower $\{\Gamma^n\}_{n\ge 0}$ converges to the identity on simply connected spaces. This shows that the layers of the tower have to contribute something to the linear part given by reduced homology $\wt{H}_*$ in order to make it up to the first derivative of the identity given by stable homotopy $\pi_*^{\rm st}$.
\end{remark}

\begin{corollary}\label{n-excisive 2}
The functor $\Gamma^n$ is $n$-excisive.
\end{corollary}

\begin{proof}
There is a fiber sequence
    $$ \Gamma^n\to\Gamma^{n-1}\to B\ol{\Gamma}^n$$
of functors. By induction the statement follows from \ref{n-excisive 1}.
\end{proof}

\section{Relation to ordinary nilpotent groups}

Before we talk about ordinary nilpotent groups, we need to gather some remarks on the Goodwillie tower of the identity.
\begin{remark}\label{convergence of Goodwillie tower}
The identity functor is $1$-analytic \cite{Goo:calc2}, which shows that its Goodwillie tower converges on simply connected spaces to the identity. However, on connected spaces $X$ it converges to the Bousfield-Kan completion of $X$:
\begin{equation}\label{Bousfield-Kan}
    \holim_n P_n(\id)(X)\simeq\mathbbm{Z}_{\infty}X 
\end{equation}
This is proved on the last page of \cite{Arone-Kankaan:iterated-Snaith}.
\end{remark}

\begin{definition}
Let $D_n(\id)(\bigvee_k S^1)$ be the homotopy fiber of the map
    $$ P_n(\id)(\bigvee_k S^1)\to P_{n-1}(\id)(\bigvee_k S^1) $$
in the Goodwillie tower. Then $D_n(\id)(\bigvee_k S^1)=\Omega^{\infty}E_n$ for the following spectrum
    $$ \left(\partial^n(\id)\wedge (\bigvee_k S^1)^{\wedge n}\right)_{h\Sigma_n}=:E_n, $$
where $\partial^n(\id)$ is the $n$-th derivative of the identity.
\end{definition}

\begin{lemma}\label{degree 1}
The rational homology of $E_n$ is concentrated in degree $1$. 
\end{lemma}

\begin{proof}
There is a Serre spectral sequence
   $$ H_i(\Sigma_n, H\mathbbm{Q}_j( \partial^n(\id) \wedge(\bigvee_k S^1)^{\wedge n})) \Longrightarrow H\mathbbm{Q}_{i+j}(E_n) $$
The homology of $\Sigma_n$ with coefficients in a rational vector space vanishes for $i>0$ and gives the formula:
    $$ H\mathbbm{Q}_*E_n\cong H\mathbbm{Q}_*(\partial^n(\id))\otimes_{\mathbbm{Q}[\Sigma_n]}H\mathbbm{Q}_*((\bigvee_k S^1)^{\wedge n}).$$
The homology of $(\bigvee_k S^1)^{\wedge n}$ is concentrated in degree $n$. 
By the work of Johnson \cite{Johnson:thesis} and Arone-Mahowald \cite{AroMah:id} we know that the spectrum $\partial^n(\id)$ is non-equivariantly equivalent to $\bigvee_{(n-1)!}S^{1-n}$ and so has homology concentrated in degree $1-n$. 
So $1-n+n=1$, and the statement follows.
\end{proof}

\begin{lemma}\label{pi0 = n-nilpotent}
The group $\pi_1 P_n(\id)(\bigvee_k S^1)$ is nilpotent of degree $n$.
\end{lemma}

\begin{proof}
We can settle the case $n=1$ right away:
   $$ \pi_0\Omega P_1(\id)(\bigvee_k S^1)\cong\pi_1(\Omega^{\infty}\Sigma^{\infty}\bigvee_k S^1)\cong\mathbbm{Z}^k\cong F_k/\Gamma_2F_k $$
In particular, this group is nilpotent of degree $1$. We proceed by induction on $n$.
From the Goodwillie tower we have for each $n\ge 1$ the following exact sequence of groups:
\begin{align}\begin{split}\label{exact sequence}
      \pi_1 D_n(\id)(\bigvee_k S^1) &\to\pi_1 P_n(\id)(\bigvee_k S^1)\to\pi_1 P_{n-1}(\id)(\bigvee_k S^1) \\
                                    &\to\pi_0 D_n(\id)(\bigvee_k S^1)\cong 0
\end{split}\end{align}
The last group vanishes because of \ref{degree 1}.
By Goodwillie's results \cite{Goo:calc3} the spaces $D_n(\id)(\bigvee_k S^1)$ are infinite loop spaces and the map in the Goodwillie tower is a principal fibration, i.e. there is a homotopy pullback square
\diagr{  P_n(\id)(\bigvee_k S^1)\ar[r]\ar[d] & \ast \ar[d] \\
         P_{n-1}(\id)(\bigvee_k S^1) \ar[r] & BD_n(\id)(\bigvee_k S^1) }
where $BD_n(\id)(\bigvee_k S^1)$ is a delooping of $D_n(\id)(\bigvee_k S^1)$ and therefore simply connected. Let 
    $$K_n= {\rm im}\left[\pi_1 D_n(\id)(\bigvee_k S^1)\to\pi_1 P_n(\id)(\bigvee_k S^1)\right] .$$
Then the short exact sequence
\begin{align}\begin{split}\label{exact sequence2}
      0\to K_n  &\to\pi_1 P_n(\id)(\bigvee_k S^1)\to\pi_1 P_{n-1}(\id)(\bigvee_k S^1) \to 0
\end{split}\end{align}
is a central extension. 
It follows inductively that the group $\pi_1 P_n(\id)(\bigvee_k S^1)$ is nilpotent of degree $n$. 
\end{proof}

\begin{lemma}\label{torsion}
The groups $\pi_sP_n(\id)(\bigvee_k S^1)$ are finite for $s\ge 2$.
\end{lemma}

\begin{proof}
These groups are finitely generated. 
So it is enough to prove that the groups $\pi_sP_n(\id)(\bigvee_k S^1)$ are torsion above degree $s=1$. 
We will prove this by induction along the Goodwillie tower where the case $n=0$ is obvious, because the space is contractible.
Next we know by \ref{degree 1} that the rational homology of $D_n(\id)(\bigvee_k S^1)$ is concentrated in degree $1$. Since $D_n(\id)(\bigvee_k S^1)$ is an infinite loop space, a form of the Hurewicz theorem \cite[thm. 9.6.20]{Spanier} tells us that also the groups $\pi_sD_n(\id)(\bigvee_k S^1)$ are torsion for $s> 1$. The result now follows from the long exact homotopy sequence of the Goodwillie tower. 
\end{proof}

\begin{corollary}\label{limit iso}
There is an isomorphism of groups:
   $$ \pi_1\holim_nP_n(\id)(\bigvee_k S^1)\cong\lim_n\pi_1P_n(\id)(\bigvee_k S^1) $$
\end{corollary}

\begin{proof}
There is the Milnor exact sequence:
   $$ 0\to\lim_n\!\mbox{}^1\pi_2P_n(\id)(\bigvee_kS^1)\to\pi_1\holim_nP_n(\id)(\bigvee_kS^1)\to\lim_n\pi_1P_n(\id)(\bigvee_kS^1)\to 0 $$
By \ref{torsion} the $\lim^1$-term vanishes.
\end{proof}

The following conjecture is related to the vanishing of $\lim^1\pi_2P_n(\id)(\bigvee_kS^1)$.
\begin{conjecture}[Arone-Mahowald-Kuhn] 
For each prime $p$ the map 
   $$\pi_sP_{p^n}(\id)(S^1)_{(p)}\to\pi_sP_{p^{n-1}}(\id)(S^1)_{(p)}$$ 
is null for $s\ge 2$. 
\end{conjecture}

Now we can describe the relation of $\mc{P}_n$ to the set-valued theory of ordinary $n$-nilpotent groups $\Nil_n$, whose $k$-ary operations are given by the free $n$-nilpotent group on $k$ generators:
    $$ \Nil_n(k^+,1^+)=F_k/\Gamma_{n+1}F_k $$
Here $F_k$ is the free group on $k$ generators. 
We can exhibit this theory by applying $\pi_0$ to the theory of homotopy $n$-nilpotent groups. First observe that for the case $n=\infty$ the statement
    $$ \pi_0\mc{P}_\infty(k^+,1^+)=\pi_1(\bigvee_k S^1)\cong F_k\cong\Nil_{\infty}(k^+,1^+) $$
follows from the Seifert-Van Kampen theorem. This isomorphism has an analogue for finite $n$. 
There is a map
    $$ F_k\cong\pi_1(\bigvee_k S^1)\to\pi_1P_n(\id)(\bigvee_k S^1)\cong\pi_0\mc{P}_n(k^+,1^+) $$
induced by the natural transformation $\id\to P_n(\id)$, which factors like
    $$ \alpha_n\co F_k/\Gamma_{n+1}F_k\to\pi_1P_n(\id)(\bigvee_k S^1),$$
because the target is $n$-nilpotent by \ref{pi0 = n-nilpotent}.

\begin{theorem}\label{pi0-iso}
We have an isomorphism of groups:
\diagr{  \alpha_n\co\Nil_n(k^+,1^+)=F_k/\Gamma_{n+1}F_k \ar[r]^-{\cong} & \pi_1P_n(\id)(\bigvee_k S^1)\cong\pi_0\mc{P}_n(k^+,1^+) }
This induces an isomorphism of categories $\Nil_n\cong\pi_0\mc{P}_n$.
\end{theorem}

\begin{proof}
We show first that $\alpha_n$ is injective by constructing a left inverse $\beta_n$.
According to \ref{n-excisive 2} the functor $\Gamma^n$ is $n$-excisive. So there is a natural transformation $P_n(\id)\to\Gamma^n$ under the identity functor. If we evaluate this diagram on $\bigvee_kS^1$ and apply $\pi_1\cong\pi_0G$, we obtain a map $b_n$ making the following diagram commutative:
\diagr{ F_k \ar[rr]^-{\gamma_{n+1}} \ar[dd]^{\cong} && F_k/\Gamma_{n+1}F_k \ar[dd]^-{\cong}_-{f_n} \ar[dl]^-{\alpha_n} \\
          & \pi_1P_n(\id)(\bigvee_kS^1) \ar[dr]_-{b_n} & \\ 
        \pi_0G(\bigvee_kS^1) \ar[rr]\ar[ur] & & \pi_0\left[G(\bigvee_kS^1)/\Gamma_{n+1}G(\bigvee_kS^1)\right] }
Let $\beta_n=f_n^{-1}b_n$.
It follows that $\beta_n\alpha_n\gamma_{n+1}=\gamma_{n+1}$. Since $\gamma_{n+1}$ is the universal map into an $n$-nilpotent group, we have $\beta_n\alpha_n=\id$.

But the map $\alpha_n$ is also surjective.
Let $Q_n$ be the quotient of $\alpha_n$, i.e. the pointed set of left cosets. We obtain a short exact sequence of towers:
\diagr{          & \vdots  \ar[d]                 & \vdots \ar[d]                             & \vdots \ar[d] & \\
        0 \ar[r] &  F_k/\Gamma_{n+1}F_k\ar[r]^-{\alpha_{n+1}}\ar[d] & \pi_1P_{n}(\bigvee_k S^1) \ar[r]\ar[d] & Q_{n} \ar[r]\ar[d] & 0 \\
        0 \ar[r] &  F_k/\Gamma_{n}F_k\ar[r]^-{\alpha_n}\ar[d] & \pi_1P_{n-1}(\bigvee_k S^1) \ar[r]\ar[d] & Q_{n-1} \ar[r]\ar[d] & 0 \\
                 &  \vdots                        & \vdots                                    & \vdots     &  } 
From (\ref{exact sequence}) it follows that the vertical maps
are surjective for $n\ge 1$. 
So all $\lim^1$-terms vanish. In the limit we obtain a short exact sequence:
    $$ 0\to\lim_n F_k/\Gamma_{n+1}F_k\to\lim_n\pi_1P_n(\id)(\bigvee_k S^1)\to\lim_nQ_n\to 0 $$
We can combine the weak equivalence (\ref{Bousfield-Kan}) with the isomorphism from \ref{limit iso} to conclude that $\lim_nQ_n\cong\ast$. In turn we have $Q_n=\ast$ for all $n\ge 1$, since all tower maps are surjective. So each $\alpha_n$ is an isomorphism.
\end{proof}

\section{The Goodwillie tower of the identity}

The natural coaugmentation $\id\to P_n(\id)$ induces a functor 
    $$p_n\co\mc{P}_{\infty}\to\mc{P}_n$$ 
of simplicial theories. 

\begin{definition}
We will denote the functor obtained by pulling back along the functor $p_n$ by
    $$ \varphi_n\co\mc{S}_*^{\mc{P}_n}\to\mc{S}_*^{\mc{P}_{\infty}}. $$
This is just the forgetful functor. It has an $\mc{S}_*$-left adjoint $\lambda_n\co\mc{S}_*^{\mc{P}_{\infty}}\to\mc{S}_*^{\mc{P}_{n}}$.
\end{definition}

\begin{remark}
One can easily prove that the pair $(\varphi_n,\lambda_n)$ forms a Quillen pair for the homotopy algebra model structures on both sides.
\end{remark}

\begin{definition}
We let 
    $$ L\lambda_n\co\mc{S}_*^{\mc{P}_{\infty}}\to\mc{S}_*^{\mc{P}_{n}}$$
be the enriched homotopy left Kan extension, which is obtained by precomposing $\lambda_n$ with a projective cofibrant replacement functor on $\mc{S}^{\mc{P}_{\infty}}$. 
\end{definition}

\begin{theorem}\label{main theorem}
Let $\dgrm{X}$ be a local homotopy $\mc{P}_{\infty}$-algebra with $\dgrm{X}(1^+)\simeq \Omega X$ for some reduced simplicial set $X$.
Then there is a natural weak equivalence
    $$ (L\lambda_n\dgrm{X})(1^+)\simeq \Omega P_n(\id)(X). $$
\end{theorem}

\begin{proof}
In the case $X=\bigvee_kS^1$ the associated homotopy $\mc{P}_{\infty}$-algebra $\dgrm{X}$ is given by
    $$ \dgrm{X}[\ell]= \Omega(\bigvee_k S^1)^\ell\simeq \mc{P}_{\infty}(k^+,\ell^+), $$
while 
    $$ \mc{P}_n(k^+,1^+)\simeq \Omega P_n(\id)(\bigvee_kS^1). $$ 
Enriched left Kan extension preserves representable functors, so we have an isomorphism:
    $$ L\lambda_n\mc{P}_{\infty}(k^+,\free)\cong\mc{P}_n(k^+,\free) $$
This proves the case $X=\bigvee_kS^1$.
Now we observe that every reduced finite simplicial set $X$ is weakly equivalent to the realization of a bisimplicial set $X_\bullet$, which consists degreewise of a finite wedge of copies of the circle $S^1$. The statement now follows from the next theorem \ref{realization} applied with $F=\Omega P_n(\id)$ and $r=k=1$.
\end{proof}

\begin{remark}
The previous theorem can be viewed as stating that $\Omega P_n(\id)(X)$ is the free homotopy $n$-nilpotent group on the loop space $\Omega X$.
\end{remark}

The condition $E_n(k)$ in the next theorem is defined in \cite[Def. 4.1]{Goo:calc2}. The theorem itself is taken from the unpublished version \cite{Mauer-Oats:unpub} of Mauer-Oats' thesis. The published version is \cite{Mauer-Oats}.

\begin{theorem}\label{realization}
Let $F$ be a reduced finitary analytic functor from spaces to spaces satisfying condition $E_n(rn-c)$ for all $n\ge 1$. If $X_\bullet$ is a simplicial $k$-connected space with $k\ge\max(r,-c)$ then
    $$ F|X_\bullet|\simeq|FX_\bullet|. $$
\end{theorem}

A sketch of the proof goes as follows: First one observes that homogeneous functors with connective coefficient spectrum commute with realizations. Then the theorem follows by induction up the Goodwillie tower. All along one checks that the connectivity estimates allow one to apply theorem \cite[B.4]{BF:gamma} that gives sufficient conditions for the realization functor to commute with pullbacks.

\section{Values of $n$-excisive functors}

Finally we will prove that functors of the form $\Omega F$ with $F$ $n$-excisive naturally take values in the category of homotopy $n$-nilpotent groups. We take this as a justification of the usefulness of the notion of homotopy $n$-nilpotent groups.

We need to compose functors.
But two functors $F$ and $G$ in \mc{F} cannot be composed directly. However, we can extend the functor $F\co\Sfinp\to\mc{S}_*$ to a functor $\mc{S}_*\to\mc{S}_*$ by enriched left Kan extension. By abuse of language we will denote this functor again by $F$. Then the composition $F\circ G$ is well-defined.

Observe that the functor $\free\circ G\co\mc{F}\to\mc{F}$ commutes with finite limits. And there is a functor 
   $$ P_n(\free\circ G)\co\mc{F}\to\mc{F}, F\mapsto P_n(F\circ G),$$
which also commutes with finite limits.
The map $F\to P_nF$ induces a map
\begin{equation}\label{Pn and composition1}
   P_n(F\circ G) \to P_n((P_nF)\circ G) 
\end{equation} 
under $F\circ G$.

\begin{lemma}\label{Pn and composition2}
The map \emph{(\ref{Pn and composition1})} consists of objectwise weak equivalences.
\end{lemma}

The proof is taken from Michael Ching in \cite[Prop. 6.1(1)]{Ching:spectra-chain-rule}, who formulates the statement for functors from spectra to spectra. But the proof goes through for our case of functors from pointed spaces to pointed spaces.

\begin{theorem}\label{values}
Let $F$ be a functor of the form $F=\Omega G$ with an $n$-excisive functor $G$ in \mc{F}. Then for any $X$ in $\Sfinp$ the space $F(X)$ is a homotopy $n$-nilpotent group.
\end{theorem}

\begin{proof}
Recall the full inclusion functor $I_n\co\mc{P}_n\to\mc{F}$, which obviously commutes with products. The same applies to the functor $P_n(\free\circ G)\co\mc{F}\to\mc{F}$. Observe also that both functors commute up to homotopy with $\Omega$.
We arrive at the conclusion that for any $F=\Omega G$ the functor $\mc{P}_n\to\mc{S}_*$ given by
\begin{align*}
      \prod_{k}\Omega P_n(\id) &\mapsto P_n((\prod_k\Omega P_n(\id))\circ G) \\
                               &\simeq\prod_k\Omega P_n(P_n(\id)\circ G)\simeq\prod_k\Omega P_n(G)\simeq F^k 
\end{align*}
preserves products up to weak equivalence by \ref{Pn and composition2}. For any space $X$ the evaluation functor $\Ev_X\co\mc{F}\to\mc{S}_*$ also preserves products. Hence $F(X)$ is a homotopy $n$-nilpotent group.
\end{proof}

The following theorem is formal.
\begin{theorem}\label{enriched}
Let $G$ be a functor of the form $G=\Omega H$ with an $n$-excisive functor $H$ in \mc{F}. Let $F$ be an arbitrary functor in \mc{F}. Then $\mc{F}(F,G)$ is a homotopy $n$-nilpotent group naturally in $F$ and $G$.
\end{theorem}


\begin{thebibliography}{10}

\bibitem{Arone-Kankaan:iterated-Snaith}
G.~Arone and M.~Kankaanrinta.
\newblock A functorial model for iterated {S}naith splitting with applications
  to calculus of functors.
\newblock In {\em Stable and unstable homotopy (Toronto, ON, 1996)}, volume~19
  of {\em Fields Inst. Commun.}, pages 1--30. Amer. Math. Soc., Providence, RI,
  1998.

\bibitem{AroMah:id}
G.~Arone and M.~Mahowald.
\newblock The {G}oodwillie tower of the identity functor and the unstable
  periodic homotopy of spheres.
\newblock {\em Invent. Math.}, 135(3):743--788, 1999.

\bibitem{Badzioch:theories}
B.~Badzioch.
\newblock Algebraic theories in homotopy theory.
\newblock {\em Ann. of Math. (2)}, 155(3):895--913, 2002.

\bibitem{Badzioch-Chung-Voronov}
B.~Badzioch, K.~Chung, and A.~A. Voronov.
\newblock The canonical delooping machine.
\newblock {\em J. Pure Appl. Algebra}, 208(2):531--540, 2007.

\bibitem{BCR:calc}
G.~Biedermann, B.~Chorny, and O.~R\"ondigs.
\newblock {Goodwillie's Calculus and model categories}.
\newblock {\em Advances in Mathematics}, 214:92--115, 2007.

\bibitem{BF:gamma}
Bousfield and Friedlander.
\newblock {Homotopy theory of $\Gamma$-spaces, spectra, and bisimplicial sets}.
\newblock In {\em {Geometric Applications of Homotopy Theory II}}, number 658
  in {Lecture Notes in Mathematics}. Springer, 1978.

\bibitem{Ching:spectra-chain-rule}
M.~Ching.
\newblock {A chain rule for Goodwillie derivatives of functors from spectra to
  spectra}.
\newblock arXiv:math.AT/07105567, 2007.

\bibitem{Curtis:relations}
E.~Curtis.
\newblock Some relations between homotopy and homology.
\newblock {\em {The Annals of Mathematics}}, 82(3):386--413, 1965.

\bibitem{Dold:free}
A.~Dold.
\newblock Homology of symmetric products and other functors of complexes.
\newblock {\em Ann. of Math. (2)}, 68:54--80, 1958.

\bibitem{Dwyer:localizations}
W.~G. Dwyer.
\newblock Localizations.
\newblock In {\em Axiomatic, enriched and motivic homotopy theory}, volume 131
  of {\em NATO Sci. Ser. II Math. Phys. Chem.}, pages 3--28. Kluwer Acad.
  Publ., Dordrecht, 2004.

\bibitem{DK:equivalence}
W.~G. Dwyer and D.~M. Kan.
\newblock Equivalences between homotopy theories of diagrams.
\newblock In {\em Algebraic topology and algebraic $K$-theory (Princeton, N.J.,
  1983)}, volume 113 of {\em Ann. of Math. Stud.}, pages 180--205. Princeton
  Univ. Press, Princeton, NJ, 1987.

\bibitem{GoJar:simp}
P.~G. Goerss and J.~F. Jardine.
\newblock {\em Simplicial homotopy theory}, volume 174 of {\em Progress in
  Mathematics}.
\newblock Birkh\"auser Verlag, Basel, 1999.

\bibitem{Goo:calc2}
T.~G. Goodwillie.
\newblock Calculus. {II}. {A}nalytic functors.
\newblock {\em $K$-Theory}, 5(4):295--332, 1991/92.

\bibitem{Goo:calc3}
T.~G. Goodwillie.
\newblock Calculus. {III}. {T}aylor series.
\newblock {\em Geom. Topol.}, 7:645--711 (electronic), 2003.

\bibitem{Hall:higher-commutators}
M.~Hall, Jr.
\newblock A basis for free {L}ie rings and higher commutators in free groups.
\newblock {\em Proc. Amer. Math. Soc.}, 1:575--581, 1950.

\bibitem{Hir:loc}
P.~S. Hirschhorn.
\newblock {\em Model categories and their localizations}, volume~99 of {\em
  Mathematical Surveys and Monographs}.
\newblock American Mathematical Society, Providence, RI, 2003.

\bibitem{jar:simplicial-presheaves}
J.~F. Jardine.
\newblock Simplicial presheaves.
\newblock {\em J. Pure Appl. Algebra}, 47(1):35--87, 1987.

\bibitem{Johnson:thesis}
B.~Johnson.
\newblock The derivatives of homotopy theory.
\newblock {\em Trans. Amer. Math. Soc.}, 347(4):1295--1321, 1995.

\bibitem{Joyal:letter}
A.~Joyal.
\newblock {Letter to A. Grothendieck}.

\bibitem{Kuhn:overview}
N.~J. Kuhn.
\newblock {Goodwillie towers and chromatic homotopy: an overview}.

\bibitem{Lawvere:semantics}
F.~Lawvere.
\newblock Functorial semantics of algebraic theories.
\newblock {\em Proc. Nat. Acad. Sci. U.S.A.}, 50:869--872, 1963.

\bibitem{Mauer-Oats:unpub}
A.~Mauer-Oats.
\newblock {Goodwillie Calculi}.
\newblock thesis, unpublished, UIUC, 2002.

\bibitem{Mauer-Oats}
A.~Mauer-Oats.
\newblock Algebraic goodwillie calculus and a cotriple model for the remainder.
\newblock {\em Trans. Amer. Math. Soc.}, 358(5):1869--1895, 2006.

\bibitem{Reedy:theories}
C.~Reedy.
\newblock Homology of algebraic theories.
\newblock Ph.D. Thesis, University of California, San Diego, 1974.

\bibitem{Schwede:theories}
S.~Schwede.
\newblock Stable homotopy of algebraic theories.
\newblock {\em Topology}, 40:1--41, 2001.

\bibitem{Serre:Lie}
J.-P. Serre.
\newblock {\em {Lie algebras and Lie groups}}.
\newblock Lecture Notes in Mathematics, No. 1500. Springer-Verlag, Berlin,
  2006.
\newblock 1964 lectures given at Harvard University. Corrected fifth printing
  of the second (1992) edition.

\bibitem{Spanier}
E.~H. Spanier.
\newblock {\em Algebraic topology}.
\newblock McGraw-Hill Book Co., New York, 1966.

\end{thebibliography}

\end{document}